\documentclass{article}
\usepackage[utf8]{inputenc}
\usepackage[english]{babel}
\usepackage{mathrsfs,fullpage}
\usepackage{amsthm,amssymb,amsfonts,amsmath}
\usepackage{authblk}
\usepackage{hyperref}
\usepackage{color}
\usepackage[letterpaper, portrait, margin=1in]{geometry}
\usepackage{hyperref}
\usepackage{verbatim,varioref}
\usepackage{amsmath}
\usepackage{longtable}
\usepackage{booktabs}
\usepackage{multirow}
\usepackage[table]{xcolor}
\usepackage{amsfonts}
\usepackage{listings}
\usepackage{color}
\usepackage{float}
\usepackage{tikz}
\usepackage{hyperref}
\usetikzlibrary{shapes,snakes}

\usepackage{tkz-graph}
\tikzstyle{vertex}=[circle, draw, inner sep=0pt, minimum size=6pt]
\newcommand{\vertex}{\node[vertex]}
\tikzstyle{vtx}=[circle, draw, inner sep=0pt, minimum size=12pt]
\newcommand{\vtx}{\node[vtx]}
\definecolor{darkgreen}{cmyk}{.9,0,.9,.2}

\newtheorem{theorem}{Theorem}[section]

\newtheorem{proposition}[theorem]{Proposition}

\theoremstyle{definition}
\newtheorem{definition}{Definition}
\definecolor{lightgray}{gray}{0.7}
\definecolor{midgray}{gray}{.9}
\definecolor{dkgreen}{rgb}{0,0.6,0}
\definecolor{gray}{rgb}{0.5,0.5,0.5}
\definecolor{mauve}{rgb}{0.58,0,0.82}
\definecolor{lightgray}{gray}{0.7}
\definecolor{midgray}{gray}{.9}

\definecolor{lightgray}{gray}{9}

\def\g{\mathfrak{g}}
\def\h{\mathfrak{h}}
\def\sp{\mathfrak{sp}_4(\mathbb{C})}
\def\sl{\mathfrak{sl}_3(\mathbb{C})}
\def\hroot{\tilde\alpha}
\def\a{\alpha}
\def\w{\varpi}

\def\NN{{\mathbb N}}
\def\ZZ{{\mathbb Z}}
\def\QQ{{\mathbb Q}}
\def\RR{{\mathbb R}}
\def\FF{{\mathbb F}}
\def\CC{{\mathbb C}}
\newcommand{\cupdot}{\mathbin{\mathaccent\cdot\cup}}
\def\hroot{\tilde{\alpha}}
\def\B{\mathcal{B}}
\def\C{\mathcal{C}}
\def\hr{\tilde{\alpha}}
\def\A{\mathcal{A}}
\def\D{\mathcal{D}}
\def\a{\alpha}
\let\oldlongtable\longtable
\let\endoldlongtable\endlongtable
\renewenvironment{longtable}{\rowcolors{2}{white}{midgray}\oldlongtable} {
\endoldlongtable}
\definecolor{lightgray}{gray}{0.55}
\title{Computing weight $q$-multiplicities for the representations of the simple Lie algebras} 
 \author[1]{Pamela E. Harris\thanks{Corresponding author: \textcolor{blue}{\href{mailto:peh2@williams.edu}{peh2@williams.edu}}. The first and third author were supported by NSF award DMS-1620202.}}
 \author[2]{Erik Insko\thanks{\textcolor{blue}{\href{mailto:einsko@fgcu.edu}{einsko@fgcu.edu}}}} 
 \author[3]{Anthony Simpson\thanks{\textcolor{blue}{\href{mailto:als7@williams.edu}{als7@williams.edu}}}}
 \affil[1,3]{Department of Mathematics and Statistics, Williams College}
 \affil[2]{Department of Mathematics, Florida Gulf Coast University}
 \renewcommand\Authands{ and }
\date{}
\begin{document}
\maketitle

\abstract{
The multiplicity of a weight $\mu$ in an irreducible representation of a simple Lie algebra $\g$ with highest weight $\lambda$ can be computed via the use of Kostant's weight multiplicity formula. This formula is an alternating sum over the Weyl group and involves the computation of a partition function. In this paper we consider a $q$-analog of Kostant's weight multiplicity and present a SageMath program to compute  $q$-multiplicities for the simple Lie algebras. }\\


\section{Introduction}

Throughout this paper we let $\mathfrak g$ be a simple Lie algebra of rank $r$ and we let $\mathfrak h$ be a Cartan subalgebra of $\mathfrak{g}$. We let $\Phi$ be the set of roots corresponding to $(\mathfrak {g,h})$, and let $\Phi^+\subseteq\Phi$ be the set of positive roots, while $\Delta\subseteq\Phi^+$ denotes the set of simple roots. We let $P(\mathfrak{g})$ be the set of integral weights, while $P_+(\mathfrak{g})$ denotes the set of dominant integral weights. Let $W$ denote the Weyl group, and for any $w\in W$, we let $\ell(w)$ denote the length of $w$. Recall that $W$ is generated by the reflections $s_1,\ldots,s_r$, where $s_i$ is the root reflection corresponding to the simple root $\alpha_i$$\in\Delta$.
For good general references see \cite{GW,Knapp}.

We let $L(\lambda)$ denote the irreducible (complex) representation of $\g$ with highest weight $\lambda$ and we let $\mu$ be a weight of this representation. To compute the $q$-multiplicity of the weight $\mu$ in $L(\lambda)$ we use the $q$-analog of Kostant's weight multiplicity, defined by Lusztig  \cite{LL}, as
\begin{align}
m_q(\lambda,\mu)&=\sum_{\sigma\in W}(-1)^{\ell(\sigma)}\wp_q(\sigma(\lambda+\rho)-\rho-\mu).\label{eq:qKWMF}
\end{align}
For any weight $\xi\in\mathfrak{h}^*$, we have that $\wp_q(\xi)=c_0+c_1q+c_2q^2+c_3q^3+\cdots+c_kq^k,$ where $c_i$ is the number of ways to 
write $\xi$ as a nonnegative integral linear combination of exactly $i$ positive roots. The $q$-analog of Kostant's weight multiplicity formula when evaluated at $q=1$ recovers the weight multiplicity formula \cite{KMF}, which is given by
\begin{align}
m(\lambda,\mu)=\displaystyle\sum_{\sigma\in W}^{}(-1)^{\ell(\sigma)}\wp(\sigma(\lambda+\rho)-(\mu+\rho))\label{eq:KWMF},
\end{align} where $\wp$ denotes Kostant's partition function and $\rho=\frac{1}{2}\sum_{\alpha\in\Phi^+}\alpha$. Kostant's partition function is the nonnegative integer-valued function, $\wp$, 
defined on $\mathfrak h^*$, by $\wp(\xi)$ = number of ways
$\xi$ may be written as a nonnegative integral linear combination of positive roots, for $\xi\in\mathfrak{h}^*$. Hence it is clear that $\wp_q|_{q=1}=\wp$ and thus $m_q(\lambda,\mu)|_{q=1}=m(\lambda,\mu)$.

Although a formula to compute $q$-multiplicities exists, its implementation is often difficult. The main obstructions include the fact that for a Lie algebra of rank $r$ the order of the Weyl group, indexing the number of terms in Equation \eqref{eq:qKWMF}, is factorial in $r$. Additionally the partition function involved has no known closed formulas, albeit for some very special cases \cite{Harris1,HIO}. Even with these complications there has been much work done in computing weight multiplicities \cite{BBCV,BV,BGR,Deckhart,FGP,Harris3,HarrisThesis,Harris2, Harris1,SB} and algorithms associated to this computation \cite{Barvinok,Barvinok2,BP,BZ,Cochet}. However, less is known about the $q$-analog computation, with the following work in that direction \cite{Gupta}.

This paper presents a SageMath program to compute the $q$-multiplicities of weights in highest weight representations of the simple Lie algebras. The source code for the program is found in Appendix \ref{sec:code} and the program can be downloaded from GitHub \textcolor{blue}{\href{https://github.com/antman1935/lie_algebras}{here}}.
 Our program computes this multiplicity by exploiting the observation that, in practice, most terms in \eqref{eq:qKWMF} and \eqref{eq:KWMF} are zero. Hence, we reduce the weight multiplicity computation by determining the Weyl group elements that contribute nontrivially to these equations. With this in mind, we give the following definition.

\begin{definition}\label{WAS}
 For $\lambda,\mu$ dominant integral weights of $\mathfrak g$ define the \emph{Weyl alternation set} to be \[\mathcal A(\lambda,\mu)=\{\sigma\in W|\;\wp(\sigma(\lambda+\rho)-(\mu+\rho))>0\}.\]
\end{definition}

There has been recent work in computing the Weyl alternation sets for certain special weights. For example, Harris showed that the number of Weyl group elements contributing nontrivially to the multiplicity of the zero weight in the adjoint representation (the representation with highest weight equal to the highest root) of $\mathfrak{sl}_{r+1}(\mathbb{C})$ was given by Fibonacci numbers \cite{Harris1}. Later, Harris, Insko, and Williams, generalized these results to show that the number of Weyl group elements contributing nontrivially to the multiplicity of the zero weight in the adjoint representation of all classical Lie algebras is governed by linear homogeneous recurrence relations with constants coefficients \cite{HIW}. 

As an application of our program we computed the Weyl alternation sets corresponding to the Weyl group elements contributing nontrivially to the multiplicity of the zero weight in the adjoint representation of the exceptional Lie algebras. We then used these Weyl alternation sets to verify the following result of Lusztig for the exceptional Lie algebras \cite{LL}: If $\mathfrak{g}$ is a simple Lie algebra with highest root $\hroot$, then $m_q(\hroot,0)=\sum_{i=1}^r q^{e_i}$, where $e_1,\ldots,e_r$ are the exponents of $\mathfrak{g}$. We recall that the exponents of $\mathfrak{g}$ are related to the degrees of the basic invariants. That is, the degrees are obtained by simply increasing the exponents by one \cite{JH}. For each exceptional Lie algebra $\g$, Table \vref{tab1}, gives the order of the Weyl group, the cardinality of the Weyl alternation set $\A(\hroot,0)$ where $\hroot$ is the highest root of $\g$, and the exponents of $\g$. From this table, it is evident that computing the Weyl alternation set reduces the computation involved in $m_q(\hr,0)$ drastically. Tables of these computations can be found in Appendix \ref{sec:application}.

\begin{table}[h!]
\centering
\begin{tabular}{|c|c|c|c|}
\hline
Lie algebra&Exponents&$|W|$&$|\A(\hr,0)|$\\ \hline
$G_2$&1,5&12&2\\\hline
$F_4$&1,5,7,11&1152&25\\\hline
$E_6$&1,4,5,7,8,11&25920 &58\\\hline
$E_7$&1,5,7,9,11,13,17&2903040&258\\\hline
$E_8$&1,7,11,13,17,19,23,29&696729600&2318\\\hline
\end{tabular}
\caption{The exceptional Lie algebras' exponents, Weyl group order, and cardinality of $\A(\hroot,0)$.}
\label{tab1}
\end{table}

Another component of our program is the ability to calculate the values of $\wp({\xi)}$ and $\wp_q({\xi})$ for individual weights $\xi$ in two ways. The program can implement the partition via a recursive algorithm or via a geometric series expansion, a generalization of Euler's generating function for integer partitions. This allows the program to be a tool for those interested in partition theory, as the program can be be used for computing the value of the partition function without having to compute weight multiplicities.

\section{Computer implementation}\label{comp}
The following sections provide our computer implementation of \eqref{eq:qKWMF} and the structures used, including both versions for the computation of Kostant's partition function and its $q$-analog.

\subsection{Structures}

Using SageMath, we automated Kostant's weight multiplicity formula and the $q$-analog of Kostant's weight multiplicity formula. SageMath provides robust tools for instantiating Lie algebras, in particular, the exceptional Lie algebras, which are the object of our study. This allowed us to compute the result of applying $\sigma\in W$ to any weight in an exceptional Lie algebra. In this process, we wrote supporting code for computing the value of the partition function, and its $q$-analog. To find values of these functions, we created two classes in Python. The first, called the $\emph{Weight class}$, provides a computer representation of a weight that affords us the ability, through basic operations and condition testing, to create a structure that finds the partition of a weight as a nonnegative integral sum of positive roots. The second class (structure), called a $\emph{Partition Tree}$, computes \eqref{eq:qKWMF} and \eqref{eq:KWMF} by instantiating every possible nonnegative integral combination of positive roots equal to the weight passed as an argument to the formulas. When the tree is instantiated, its input is the weight we are trying to partition, and a list of the positive roots that can be used in the partitioning.
     
To begin the process, the root of the tree is created, then its children are instantiated. This is done by using the first item in the list of positive roots being used to partition the weight stored at the root of the tree. Note that root has two uses, the root of a tree, or a positive root. To avoid confusion, we will always refer to the root of the tree in full. If we let $\xi$ represent the weight we are trying to partition, and ${r_1}$ denote the root that we are currently using to do so, then the root has a child for every element of $\{\xi - nr_1 \mid n \in \mathbb{Z}, n \geq 0 \}$ that contains no simple roots with negative coefficients. Each child stores a unique new weight $\xi - nr_1$, and $n$, the number of times ${r_1}$ was used in the partitioning. This represents the subtraction of the first positive root from the original weight as many times as possible. From here, the process begins at each of the root's children, but using ${r_2}$, the second root in the list of positive roots. The same process starts again at each of the children of the children of the root, but with ${r_3}$, the third root in the list of positive roots. In general, a node in the tree at a depth of $i$ from the root, will have its children created through the same process using ${r_{i+1}}$, the ${(i+1)^{th}}$ root in the list of positive roots.
    
The process is complete when either the weight is reduced to the zero weight, or any one of the coefficients in the linear combination of simple roots, equal to the weight stored in the children of the tree, is negative. When it is reduced to zero, we have found a successful partitioning, and so the terminal node, a constant that is instantiated at the beginning of the program, is made to be a child of the node where the weight was reduced to zero. Not all branches of the tree will end in a terminal node, so another function goes through all tree branches and removes those that do not end in the terminal node, since they do not represent successful partitions, i.e. those weights where a coefficient becomes negative.

The following image represents the partition tree of the weight in \eqref{eq:result}. Note that we provide the list of positive roots to the right of the tree, and underline the one being used to partition at each level of the tree. Note that when we talk about the relationship between nodes in the tree, we say that if one node has a connection to another node, then those nodes have an edge between them, where the first node is the source node and the second is the destination node, with respect to that edge. In this specific tree structure, we have that the weight stored in the destination node is equal to the weight stored in the source node minus $n$ times the highlighted positive root. Note that in the description of the code, the source node holds on to the value of $n$, but for the sake of the graphic, it is easier to have them associated to the edges between nodes.

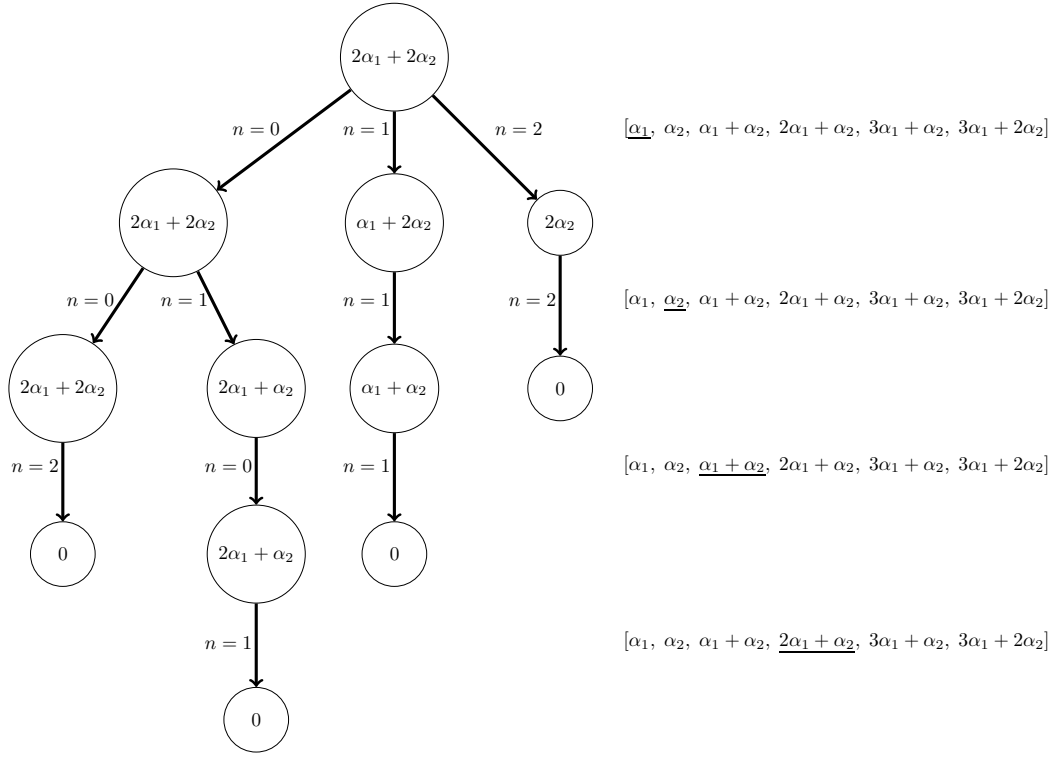
\begin{figure}[h!]
\centering
\resizebox{5.5in}{!}{%
\begin{tikzpicture}
\vertex[](a0) at (0,0) {$\;\;2\alpha_1+2\alpha_2\;\;$};
\vertex[](b1) at (-4,-3) {$\;\;2\alpha_1+2\alpha_2\;\;$};
\vertex[](b2) at (0,-3) {$\;\;\alpha_1+2\alpha_2\;\;$};
\vertex[](b3) at (3,-3) {$\;\;\;2\alpha_2\;\;\;$};
\vertex[](c1) at (-6,-6) {$\;\;2\alpha_1+2\alpha_2\;\;$};
\vertex[](c2) at (-2.5,-6) {$\;\;2\alpha_1+\alpha_2\;\;$};
\vertex[](c3) at (0,-6) {$\;\;\alpha_1+\alpha_2\;\;$};
\vertex[](c4) at (3,-6) {$\;\;\;\;\;0\;\;\;\;\;$};
\vertex[](d1) at (-6,-9) {$\;\;\;\;\;0\;\;\;\;\;$};
\vertex[](d2) at (-2.5,-9) {$\;\;2\alpha_1+\alpha_2\;\;$};
\vertex[](d3) at (0,-9) {$\;\;\;\;\;0\;\;\;\;\;$};
\vertex[](e1) at (-2.5,-12) {$\;\;\;\;\;0\;\;\;\;\;$};
\draw[ultra thick,->] (a0)--(b1);
\draw[ultra thick,->] (a0)--(b2);
\draw[ultra thick,->] (a0)--(b3);
\draw[ultra thick,->] (b1)--(c1);
\draw[ultra thick,->] (b1)--(c2);
\draw[ultra thick,->] (b2)--(c3);
\draw[ultra thick,->] (b3)--(c4);
\draw[ultra thick,->] (c1)--(d1);
\draw[ultra thick,->] (c2)--(d2);
\draw[ultra thick,->] (c3)--(d3);
\draw[ultra thick,->] (d2)--(e1);
\node at (-2.5,-1.3) {{$n=0$}} ;
\node at (-.5,-1.3) {{$n=1$}} ;
\node at (2.25,-1.3) {{$n=2$}} ;
\node at (8,-1.3) {$[\underline{\alpha_1},\;\alpha_2,\;\alpha_1+\alpha_2,\;2\alpha_1+\alpha_2,\;3\alpha_1+\alpha_2,\;3\alpha_1+2\alpha_2]$};
\node at (-5.5,-4.4) {{$n=0$}} ;
\node at (-3.8,-4.4) {{$n=1$}} ;
\node at (-.5,-4.4) {{$n=1$}} ;
\node at (2.5,-4.4) {{$n=2$}} ;
\node at (8,-4.4) {$[\alpha_1,\;\underline{\alpha_2},\;\alpha_1+\alpha_2,\;2\alpha_1+\alpha_2,\;3\alpha_1+\alpha_2,\;3\alpha_1+2\alpha_2]$};
\node at (-6.5,-7.4) {{$n=2$}} ;
\node at (-3,-7.4) {{$n=0$}} ;
\node at (-.5,-7.4) {{$n=1$}} ;
\node at (8,-7.4) {$[\alpha_1,\;\alpha_2,\;\underline{\alpha_1+\alpha_2},\;2\alpha_1+\alpha_2,\;3\alpha_1+\alpha_2,\;3\alpha_1+2\alpha_2]$};
\node at (-3,-10.6) {{$n=1$}} ;
\node at (8,-10.6) {$[\alpha_1,\;\alpha_2,\;\alpha_1+\alpha_2,\;\underline{2\alpha_1+\alpha_2},\;3\alpha_1+\alpha_2,\;3\alpha_1+2\alpha_2]$};
\end{tikzpicture}}
\label{pic}\caption{Partition tree for $2\alpha_1 + 2\alpha_2$ using the positive roots of the Lie algebra $G_2$.}\label{fig:tree}
\end{figure}

Recall that when a branch represents a successful partition, the branch ends in the terminal node. For the purpose of simplifying the image in Figure~\vref{fig:tree}, only the successful branches are provided, i.e. those branches ending in zero nodes. Using Figure~\vref{fig:tree}, we can arrive at the same partitions of $2\alpha_1 + 2\alpha_2$ as we did in Section \ref{G2ex}. Note that, if we start from the root of the tree, and follow any path to a node containing zero, we will find a partition of $2\alpha_1 + 2\alpha_2$ using the positive roots of $G_2$. The collection of all paths from the node with $2\alpha_1 + 2\alpha_2$ to a node containing zero yields all partitions of this weight. For example, the path from the root to the leftmost zero, indicates that $\alpha_1$ and $\alpha_2$ were not used, and $\alpha_1 + \alpha_2$ was used twice. Thereby accounting for $2(\alpha_1+\alpha_2)$ as a partition of $2\alpha_1 + 2\alpha_2$.
    
\subsection{Computations using the structures}

Once the partition tree has been created, we can compute $\wp(\xi)$ and $\wp_q(\xi)$. That is, $\wp(\xi)$ can be determined by counting the number of branches of the partitioning tree that end with the terminal node. Because we provided a function that removed the branches that do not represent a successful partition, every branch ends with the terminal node, so computing $\wp(\xi)$ only requires us to count the number of branches in the tree, where a branch is just any path from the root of the tree to any leaf in the tree. Note that there is only one leaf in the tree, namely the terminal node, so we need only check for that to know the end of the tree.
    
The computation for $\wp_q(\xi)$ requires that we keep track of the number of roots used in the partition. However, this was already done by storing $n$, from the expression ${\xi-nr_i}$. Because the number of roots used at each step of the partitioning is stored for every partition, simply summing the values of $n$ stored in every node of a branch gives the total number of roots used in any given partition. Moreover, counting the roots used in a partition is equivalent to determining the length of the branch when one considers the $n$ value of a child to be the distance between that child node and its parent node. Once the number of roots used is known for a single partition of the weight, the corresponding coefficient in our running calculation of $\wp_q(\xi)$ must be incremented. This means that every $c_i$ in $\wp_q(\xi) = c_0 + c_1q^1 + \cdots c_kq^k$ equals the total number of branches of length $i$ in the partitioning tree. From all branches in the tree, we can obtain a value for $\wp_q(\xi)$.
    
With the ability to obtain $\wp(\sigma(\lambda+\rho)-(\rho+\mu))$ and $\wp_q(\sigma(\lambda+\rho)-(\rho+\mu))$ for any element $\sigma\in W$, finding the values of $m(\lambda, \mu)$ and $m_q(\lambda, \mu)$ just requires us to iterate over all the elements of the Weyl group. Then, using the length of each element to determine the sign of each term's contribution to the sum in Kostant's weight multiplicity formula, yields the desired result.

\subsubsection{The Lie algebra of type $G_2$}\label{G2ex}

We consider the exceptional Lie algebra $G_2$ to illustrate our computations. The set of positive roots of the exceptional Lie algebra $G_2$ is 
\begin{equation}\label{g2phi}
\Phi^+ = \{\alpha_1, \alpha_2, \alpha_1+\alpha_2, 2\alpha_1+\alpha_2, 3\alpha_1+\alpha_2, 3\alpha_1+2\alpha_2\}
\end{equation}
where $3\alpha_1 + 2\alpha_1$ is the highest root, and $\rho = \frac{1}{2}\sum_{\alpha \in \Phi^+}\alpha=5\alpha_1+3\alpha_2$. The Weyl group of $G_2$ has 12 elements and is isomorphic to the dihedral group $D_6$. As every element in $W$ is generated by $s_1$ and $s_2$, the simple reflections associated with the simple roots $\alpha_1$ and $\alpha_2$, we can understand how any reflection acts on any weight by composing the the simple reflections and applying them to each of the simple roots. In this case we observe that
\begin{align} 
    \label{eq:s1ona1} s_1&: \alpha_1 \mapsto -\alpha_1 \\
    \label{eq:s1ona2} s_1&: \alpha_2 \mapsto 3\alpha_1+\alpha_2 \\
    \label{eq:s2ona1} s_2&: \alpha_1 \mapsto \alpha_1 + \alpha_2 \\
    \label{eq:s2ona2} s_2&: \alpha_2 \mapsto -\alpha_2.
\end{align}

To compute $s_1(\widetilde{\alpha} + \rho) - \rho$, we note that $\widetilde{\alpha} + \rho = 8\alpha_1 + 5\alpha_2$ and by \eqref{eq:s1ona1} and \eqref{eq:s1ona2}, it follows that
\begin{align}
    s_1(\widetilde{\alpha} + \rho) - \rho &= s_1(8\alpha_1+5\alpha_2) - (5\alpha_1 + 3\alpha_2)
= 8s_1(\alpha_1)+5s_1(\alpha_2) - (5\alpha_1 + 3\alpha_2)
= 2\alpha_1 + 2\alpha_2.\label{eq:result}
\end{align}

We can now apply Kostant's partition function to $2\alpha_1 + 2\alpha_2$, and by using the positive roots of the Lie algebra of type $G_2$ as given in \eqref{g2phi}, we see that one can write $2\alpha_1 + 2\alpha_2$ as a nonnegative integral sum of positive  roots in the following four ways
$2(\alpha_1) + 2(\alpha_2)$,
$2(\alpha_1 + \alpha_2)$, 
$1(\alpha_1) + 1(\alpha_2) + 1(\alpha_1 + \alpha_2)$, and $1(\alpha_1) + 1(2\alpha_1 + \alpha_2)$,
which use 4, 2, 3, and 2 positive roots, respectively. Thus $\wp_q(2\alpha_1 + 2\alpha_2) = 2q^2 + q^3 + q^4.$ To determine the value that $s_1$ contributes to \eqref{eq:qKWMF}, we note that $(-1)^{\ell(s_1)}=(-1)^1=-1$, hence $$(-1)^{\ell(s_1)}\wp(s_1(\widetilde{\alpha} + \rho) - \rho) = -2q^2 - q^3 - q^4.$$

Repeating this procedure for every remaining element of $W$, we arrive at the Table \vref{tb:G2}, which only lists the three elements of $W$ contributing nontrivially to $m_q(\widetilde{\alpha}, 0)$. From this one can verify that $m_q(\widetilde{\alpha}, 0)=q+q^5$ as expected.

\begin{table}[h]
\centering
\rowcolors{1}{}{lightgray}
\begin{longtable}{|c|c|p{4cm}|p{6cm}|}

\rowcolor{lightgray} $\sigma\in W$ \ &$\ell(\sigma)$ \ &$ \sigma(\widetilde{\alpha} + \rho) - \rho$ \ &$\wp_q(\xi)$\\

 $1$ & 0 & $ 3\alpha_{1} + 2\alpha_{2} $ & $ q^{1} + 2q^{2} + 2q^{3} + q^{4} + q^{5} $\\\hline
 $s_1$ & 1 & $ 2\alpha_{1} + 2\alpha_{2} $ & $ 2q^{2} + q^{3} + q^{4} $\\\hline
 $s_2$ & 1 & $ 3\alpha_{1} $ & $ q^{3} $\\\hline

\rowcolor{white}\multicolumn{4}{|c|}{$m_q(\widetilde{\alpha}, 0) = q + q^5$}\\

\hline
\end{longtable}
\caption{Weyl alternation sets and $q$-analog values for $G_2$.}\label{tb:G2}
\end{table}

\subsection{Alternate approach to computing Kostant's partition function} \label{AnAp}
It turns out that, while the exhaustive method detailed above works for $G_2$, $F_4$, $E_6$, and $E_7$, the approach required too much memory for $E_8$. Even with a machine that had 32 gigabytes of memory, partitioning a single weight took too much memory. This is an astounding fact considering that the computation for $E_7$ completed without a problem. However in the $E_8$ case, the partitions do not complete as we are using 120 positive roots. If we consider the number of nodes that would appear in the partition tree of a weight just when we consider either including, or not including each positive root, then we have $2^{120} \approx 1.33\times10^{36}$ nodes in the tree already. Hence, we  implement a well-known method to determine the value of the partition function, which is typically used for partitioning integers, and we adapt it for our purposes.

The new partition function method focuses on expanding a geometric series. We adapted Euler's formula for finding the number of partitions of an integer, denoted $p(n)$. We begin by showing how to derive this formula first by following the work presented in \cite{textforgenfuncs}. Begin by examining
\begin{equation}\label{exp1}
\sum_{n \geq 0}p(n)x^n,
\end{equation}
which is a formal power series, where the coefficient of $x^n$ counts the number of partitions of $n$. We will derive a formula for the generating function given in \eqref{exp1}. Define $S$ to be the set of partitions of an arbitrary integer. We denote the partition of an integer as $(1^{m_1}, 2^{m_2},\cdots, n^{m_n})$, where $m_i$ is the number of times $i$ is used in the partition. From this, we can see that for every partition, $\kappa$,
\begin{align}\label{exp2}
\kappa = &(1^0 \in \kappa\ or\ 1^1 \in \kappa\ or\ 1^2 \in \kappa\ or\ \cdots)\ and \\ \nonumber
&(2^0 \in \kappa\ or\ 2^1 \in \kappa\ or\ 2^2 \in \kappa\ or\ \cdots)\ and \\ \nonumber
&(3^0 \in \kappa\ or\ 3^1 \in \kappa\ or\ 3^2 \in \kappa\ or\ \cdots)\ and\ \cdots \nonumber.
\end{align}
From the definition of an element in $S$, the generating function for the elements of $S$ will be
\begin{equation}\label{exp3}
(x^0 + x^1 + x^{1 + 1} +\cdots)(x^0 + x^2 + x^{2 + 2} + \cdots)(x^0 + x^3 + x^{3+3})\cdots.
\end{equation}

Notice that the terms in the product of \eqref{exp3} are geometric sums and $\sum_{n=0}^{\infty}r^n = \frac{1}{1-r}$. Hence, we can rewrite \eqref{exp1} as
\begin{equation}
\sum_{n \geq 0}p(n)x^n = \frac{1}{1 - x}\cdot\frac{1}{1 - x^2}\cdot\frac{1}{1-x^3}\cdots,
\end{equation}
which can be rewritten as
\begin{equation}\label{exp4}
\sum_{n \geq 0}p(n)x^n = \lim_{n\to\infty}\prod_{i = 1}^{n}\frac{1}{1 - x^i}.
\end{equation}
Though \eqref{exp4} involves an infinite product, when we want to find the number of partitions for a specific integer, we can change the limit in \eqref{exp4} to be set to that integer. This is because no integer greater than the one that we are partitioning could be used. For example, if we look at partitions of $2$, we see that its partitions still have the general form given in \eqref{exp2}, but since no number greater than three could be used, we know that the every partition of 2 will be of the form $(1^{m_1}, 2^{m_2},1, \cdots, 1)$. This means that every term past $2$ in \eqref{exp4}, on either side of the equals sign, can be ignored.

We now adapt this approach to count the number of partitions for a weight using a certain set of positive roots. To be able to do this, we must first provide a translation from linear combinations of simple roots to products of the variables used in the series. We take a linear combination of simple roots, $a_1\alpha_1+\cdots+a_n\alpha_n$, and change it to $A_1^{a_1}\cdots A_n^{a_n}$. We illustrate the technique with an example. If we take the positive root $3\alpha_1 + 2\alpha_2$ in $G_2$, then the term in our geometric sum will be $A_1^3A_2^2$. If we denote the set of these translations of positive roots to these terms as \[\Phi^+_{var} = \{A_1^{a_1}A_2^{a_2}\cdots A_r^{a_r} : a_1\alpha_1+a_2\alpha_2+\cdots+a_r\alpha_r\in \Phi^+\},\] then the series that we are interested in expanding is
\begin{equation}\label{prod1}
\prod_{x \in \Phi^+_{var}}\frac{1}{1-qx}.
\end{equation}
We introduce the $q$ into our adaptation of \eqref{exp4}, to count the number of roots used in each partition. Note that in this product, we have geometric sums and we can expand \eqref{prod1} as
\begin{equation}\label{prod2}
\prod_{x \in \Phi^+_{var}}\sum_{n=0}^\infty (qx)^n.
\end{equation}

From here, we return to the specific example of $G_2$. In $G_2$, we know the positive roots are given by $\Phi^+ = \{\alpha_1, \alpha_2, \alpha_1 + \alpha_2, 2\alpha_1 + \alpha_2, 3\alpha_1 + \alpha_2, 3\alpha_1 + 2\alpha_2\}$, so $
\Phi^+_{var} = \{A_1, A_2, A_1A_2, A_1^2A_2, A_1^3A_2, A_1^3A_2^2\}$. We can use this to rewrite \eqref{prod2} for this specific case, which becomes
{\small
\begin{align}
&\sum_{n=0}^\infty (qA_1)^n\cdot\sum_{n=0}^\infty (qA_2)^n\cdot\sum_{n=0}^\infty (qA_1A_2)^n\cdot\sum_{n=0}^\infty (qA_1^2A_2)^n\cdot\sum_{n=0}^\infty (qA_1^3A_2)^n\cdot\sum_{n=0}^\infty (qA_1^3A_2^2)^n\label{prod4}.
\end{align}
}%
By expanding \eqref{prod4}, the value of $\wp_q(\xi)$ is the coefficient of the translation of $\xi$, just in the same way that the coefficient of $x^n$ was the number of partitions of $n$ in \eqref{exp1}. You can think of the translation as the bijection $trans(a_1\alpha_1+\cdots+a_r\alpha_r)=A_1^{a_1}\cdots A_r^{a_r}$ for any weight $a_1\alpha_1+\cdots+a_r\alpha_r$. Thus, we can limit the upper bound on the individual sums in \eqref{prod4} to the maximum number of times the associated positive root could have been used. We will again examine the weight $2\alpha_1 + 2\alpha_2$ from \eqref{eq:result}. Our translation of the weight $2\alpha_1 + 2\alpha_2$ is $A_1^2A_2^2$.  We know that $\alpha_1$, $\alpha_2$, and $\alpha_1+\alpha_2$ can only be used a maximum of 2 times, $2\alpha_1+\alpha_2$ can only be used once, and the other positive roots of $G_2$ cannot be used when partitioning $2\alpha_1+2\alpha_2$. If we are examining the weight $2\alpha_1 + 2\alpha_2$ as we did at the end of Section \ref{G2ex}, then we need only expand the polynomial
\begin{equation}\label{poly1}
(1+qA_1+q^2A_1^2)(1+qA_2+q^2A_2^2)(1+qA_1A_2+q^2A_1^2A_2^2)(1 + qA_1^2A_2)(1)(1).
\end{equation}
Notice that we do not have to expand \eqref{poly1} fully. Since we know exactly for which term we are searching, we need only multiply everything that will give us a coefficient of $A_1^2A_2^2$. Thus we can choose the terms
\begin{align}
(q^2A_1^2)(q^2A_2^2) &= q^4A_1^2A_2^2\\
(qA_1)(qA_2)(qA_1A_2) &= q^3A_1^2A_2^2\\
(qA_2)(qA_1^2A_2) &= q^2A_1^2A_2^2\\
(q^2A_1^2A_2^2) &= q^2A_1^2A_2^2.
\end{align}

The sum of these relevant terms yields $(2q^2+q^3+q^4)A_1^2A_2^2$, whose coefficient is exactly the value of $\wp_q(2\alpha_1 + 2\alpha_2)$ that we found at the end of Section \ref{G2ex}. Proceeding in this fashion with every element of the Weyl group that contributes nontrivially to Kostant's weight multiplicity formula we obtain $m_q(\widetilde{\alpha}, 0)$.

\subsection{Generating the Weyl alternation set}

Even with the faster method for computing the partition function, there were still too many elements of the Weyl group of $E_8$ for the computation to finish in a feasible amount of time. Because of this, we implemented a way to compute the alternation set $\mathcal{A}(\lambda, \mu)$ before proceeding to find the multiplicity of $\mu$ with respect to $\lambda$. This allowed us to work with a much smaller subset of the Weyl group, and thus reduce the computational time associated to the partitioning algorithms.

To compute the alternation set, we first check to see if the expression $1(\lambda + \rho) - (\rho + \mu)$ has any negative coefficients when expressed as a sum of simple roots. If so, then it cannot be partitioned using positive roots and the alternation set is the empty set. If the coefficients are all nonnegative integers, then we know that at least the identity element is in the alternation set, so we append the identity element to the alternation set. We then build the alternation set by taking the elements already in the set, concatenating by simple reflections, and computing $\sigma(\lambda + \rho) - (\rho + \mu)$, where $\sigma$ is the result of concatenating a simple reflection to an existing element of the alternation set. If the result of this computation results in a nonnegative integral combination of simple roots, then we include the element in the alternation set and iterate the process. Once we cannot append any simple reflections to any of the elements of the alternation set, we know that we have the full alternation set, as we have run through every element in the Weyl group.

This process is implemented with Python's list data structure, rather than sets, so to disambiguate the process further, we will explain in this context. We start the alternation set as a list with just the identity, then we begin iterating through the list. We will denote our index with the variable $i$, which will initially be set to 0. We start by taking the element in the alternation set at position $i$. We iterate over each simple reflection, append it to the reflection at position $i$ to make a new reflection and compute the same expression mentioned before. If the expression results in a nonnegative integral linear combination of simple roots, then the new reflection is appended to the end of the list. Notice that this element will also have simple reflections appended to it as $i$ increases and eventually reaches this element. This process is guaranteed to work because we know that once the expression $\sigma(\lambda+\rho)-(\mu+\rho)$ yields a negative coefficient when expressed as  a sum of simple roots, appending more reflections to $\sigma$ which make the length of $\sigma$ increase will only decrease the coefficients of the simple roots further, as was proven by the first and second author in \cite[Proposition 3.4]{HIW}. We present this result below, but omit its proof as it is both technical and lengthy.

\begin{proposition}\label{nobigger} Let $s_i $ denote the simple root reflection corresponding to $\alpha_i \in \Delta$. 
Write the weights  $\sigma(\rho)-\rho $, $\sigma s_i(\rho) - \rho $, and $s_i\sigma(\rho) - \rho$ as linear combinations of positive roots as follows:
\begin{itemize}
        \item $\sigma(\rho)-\rho = -  \sum_{\alpha_j \in \Delta} c_j \alpha_j$ with $c_j \geq 0$,
        \item $\sigma s_i(\rho) - \rho = -  \sum_{\alpha_j \in \Delta} d_j \alpha_j$ with $d_j \geq 0$, and
        \item $s_i\sigma(\rho) - \rho = -  \sum_{\alpha_j \in \Delta} e_j \alpha_j$ with $e_j \geq 0$.  
       \end{itemize} 
Then:
\begin{enumerate}
\item If $\ell(\sigma s_i)>\ell(\sigma)$, then  $d_j \geq c_j$ for all $j$ and $d_k> c_k$ for at least one $\alpha_k \in \Delta$.
\item If $\ell(s_i\sigma)>\ell(\sigma)$, then $e_j \geq c_j$ for all $j$ with $e_i> c_i$.
\end{enumerate}
\end{proposition}
 
Proposition \ref{nobigger} implies that by appending simple reflections that increase the length of elements known to be in the alternation set until the Weyl group element created results in a weight with at least one negative coefficient results in an exhaustive list of elements in the Weyl alternation set.

\subsection{Program optimizing}
Though this project is complete, it may be worthwhile to optimize the program. Should the partition tree program be optimized for memory usage, it could be possible to complete the $E_8$ computation. The advantage of the first program is that it allows the actual individual partitions to be listed from the partition tree. However, the method of series expansion allows for much quicker computations. Where the first method took approximately 169.33 seconds to run $E_6$, the method by series expansion only took 109.14 seconds. Notice that this was run on a build of SageMath that did not allow for the code to be run in parallel. From a computer scientist's viewpoint, it would be interesting to see the increase in speed if the first program was parallelized and possibly run on a GPU (Graphics Processing Unit), instead of a CPU, especially because of the independent nature of each branch of the partition tree.

\newpage
\appendix
\section{SageMath code}\label{sec:code}
The code described in Section \ref{comp} is provided in this appendix.

\subsection{The weight class}
\begin{lstlisting}
# Anthony Simpson, als7@williams.edu
from copy import copy

class Weight:
    def __init__(self, cs):
        # we just need an array of coefficients
        self.coefficients = cs
        
    def height(self):
        # we iterate through the coefficients while summing them
        ht = 0
        for x in self.coefficients:
            ht += x
        return ht
    
    def __sub__(self, other):
        # iterate through coefficients and compute differences piecewise
        c = copy(self.coefficients)
        for i in range(0, len(c)):
            c[i] -= other.coefficients[i]
        return Weight(c)
            
    def __add__(self, other):
        # iterate through coefficients and compute sums piecewise
        c = copy(self.coefficients);
        for i in range(0, len(c)):
            c[i] += other.coefficients[i]
        return Weight(c)
        
    def __str__(self):
        # print string form of array
        return str(self.coefficients)
    
    def equals(self, other):
        # check for equality of coefficients piecewise
        for i in range(0, len(self.coefficients)):
            if not (self.coefficients[i] == other.coefficients[i]):
                return False
        return True
    
    def __ne__(self, other):
        # definition based on equals method
        return not self.equals(other)
        
    def isZero(self):
        # check to see if coefficients are all zero
        for coeff in self.coefficients:
            if coeff != 0:
                return False
        return True
        
    def isPositive(self):
        # check to see if all coefficients are greater than or equal to zero
        for coeff in self.coefficients:
            if coeff < 0:
                return False
        return True
        
    def isNegative(self):
        # check to see if all coefficients are less than zero
        for coeff in self.coefficients:
            if coeff < 0:
                return True
        return False
    
    def hasFraction(self):
        # checks for fractions by subtracting integer part and checking for zero
        for coeff in self.coefficients:
             if ((coeff - int(coeff)) != 0):
                return True
        return False
\end{lstlisting}

\subsection{The partition tree class}
\begin{lstlisting}
# Anthony Simpson, als7@williams.edu
from Weight import Weight

TerminalNode = '0'

class PartitionTree:
    
    def __init__(self, r, rs, index, terms):
        self.weight = r;
        self.children = [];
        self.terms = terms
        
        # if all coefficients are zero, then we've found a successful partition
        if (r.isZero()):
            self.children.append(TerminalNode);
            
        else:
            if (index < len(rs)):
                # continue the partition by not including the current root
                self.children.append(PartitionTree(Weight(list(r.coefficients)), rs, index + 1, 0)); # exclude current root
                
                # start looking at combinations of current root
                newChild = r - rs[index]
                terms = 1
                
                # continue to subtract current root until a coefficient is below zero
                while (newChild.isPositive()):
                    self.children.append(PartitionTree(newChild, rs, index + 1, terms));
                    newChild = newChild - rs[index]
                    terms += 1
                    
    # deletes branches that do not end in the terminal node
    def cleanTree(self):
        if (len(self.children) == 0):
            del(self)
            return True
        
        delete = True
        for child in self.children:
            if child == TerminalNode: 
                # if a branch ends in TerminalNode, don't delete
                return False
            
            # if branch has multiple branches following, check to see if at least one ends in TerminalNode
            delete = delete and child.cleanTree()
            
        # if no branches ended in TerminalNode, then delete this node
        if (delete):
            del(self)
            return True
        else:
            return False
                    
    def countPartitions(self):
        # recursively count the number of unique partitions by counting the number of times we reach the end of a branch
        parts = 0;
        for child in self.children:
            # if the child is TerminalNode, then return 1, there are no other children
            if (child == TerminalNode):
                return 1;
            else:
                # there are multiple paths from this node, each one on a different branch
                parts += child.countPartitions();
        
        return parts
        
    def generatePq(self, e):
        # the first node need not be included in this computation, so we defer to a helper method on it's children
        for child in self.children:
            child.generatePqHelper(self.weight, 0, e)
        
    def generatePqHelper(self, prevRoot, rootsUsed, equationCoefficients):
        # sum the number of roots used in this partition
        rootsUsed += self.terms
            
        for child in self.children:
            if (TerminalNode == child):
                # update the coefficient of the correct power of q depending on the number roots used, if we are at the end of a partition
                equationCoefficients[rootsUsed] = equationCoefficients[rootsUsed] + 1
                return
            else:
                # otherwise continue to generate Pq with the children
                child.generatePqHelper(self.weight, rootsUsed, equationCoefficients)
\end{lstlisting}

\subsection{The $q$-analog of Kostant's weight multiplicity formula by exhaustion}
\begin{lstlisting}
# Anthony Simpson, als7@williams.edu
from sage.all import *
from Weight import Weight
from PartitionTree import PartitionTree
import __future__

count = []

#Constant Declarations

#SageMath's representation uses the standard basis instead of the simple roots, so here we instantiate these matrices to make a change of basis to the simple roots. I will automate the process of finding these matrices in the future to generalize this program for any Lie algebra
bChangeG2 = matrix([[2, 1, 0],
[1, 0, 0],
[1, 1, 1]])

bChangeF4 = matrix([[1, 1, 0, 0],
[2, 1, 1, 0],
[3, 1, 1, 1],
[2, 0, 0, 0]])

bChangeE6 = matrix([[0,   0,    0,    0,    0,   -2,    0,    0],
[ 1.0/2,  1.0/2,  1.0/2,  1.0/2,  1.0/2, -3.0/2,    0,    0],
[-1.0/2,  1.0/2,  1.0/2,  1.0/2,  1.0/2, -5.0/2,    0,    0],
[   0,    0,    1,    1,    1,   -3,    0,    0],
[   0,    0,    0,    1,    1,   -2,    0,    0],
[   0,    0,    0,    0,    1,   -1,    0,    0],
[   0,    0,    0,    0,    0,   -1,    1,    0],
[   0,    0,    0,    0,    0,    1,    0,    1]])

bChangeE7 = matrix([[   0,    0,    0,    0,    0,    0,   -2,    0],
[ 1.0/2,  1.0/2,  1.0/2,  1.0/2,  1.0/2,  1.0/2,   -2,    0],
[-1.0/2,  1.0/2,  1.0/2,  1.0/2,  1.0/2,  1.0/2,   -3,    0],
[   0,    0,    1,    1,    1,    1,   -4,    0],
[   0,    0,    0,    1,    1,    1,   -3,    0],
[   0,    0,    0,    0,    1,    1,   -2,    0],
[   0,    0,    0,    0,    0,    1,   -1,    0],
[   0,    0,    0,    0,    0,    0,    1,    1]])

bChangeE8 = matrix([[   0,    0,    0,    0,    0,    0,    0,    2],
[ 1.0/2,  1.0/2,  1.0/2,  1.0/2,  1.0/2,  1.0/2,  1.0/2,  5.0/2],
[-1.0/2,  1.0/2,  1.0/2,  1.0/2,  1.0/2,  1.0/2,  1.0/2,  7.0/2],
[   0,    0,    1,    1,    1,    1,    1,    5],
[   0,    0,    0,    1,    1,    1,    1,    4],
[   0,    0,    0,    0,    1,    1,    1,    3],
[   0,    0,    0,    0,    0,    1,    1,    2],
[   0,    0,    0,    0,    0,    0,    1,    1]])

#We set these dictionaries up so that the program could be written in general for any of the exceptional Lie algebras
basisChangeMatrices = {"G2": bChangeG2, "F4": bChangeF4, "E6": bChangeE6, "E7": bChangeE7, "E8": bChangeE8}
numOfRoots = {"G2": 6, "F4": 4, "E6": 36, "E7": 63, "E8": 120}
vectorDimension = {"G2": 3, "F4": 4, "E6": 8, "E7": 8, "E8": 8}
m = 0

#this eval method is used instead of the built in eval function to fix some associated rounding errors
def Eval(string):
    return eval(compile(str(string), '<string>', 'eval', __future__.division.compiler_flag))

def computeMq(name, lamb, mu):
    # initialize weyl group and vector space for the lie algebra
    lie_algebra = RootSystem(name).ambient_space()
    weyl_group = WeylGroup(name, prefix = "s")
    m = 0
    
    # equation is a vector that will hold the coefficients of Mq
    equation = vector([0 for i in range(0, numOfRoots[name] + 25)])
    
    # used to change the basis from the standard basis of R^n to simple roots
    changeBasis = basisChangeMatrices[name]
    
    # this collects all of the positive roots for use in partitioning
    poset = [vector(list(Eval(x))) for x in RootSystem(name).ambient_space().positive_roots()]
    
    # take those positive roots and change them to a linear combo of simple roots
    alphaPoset = [Weight(list(changeBasis * x)) for x in poset]
    
    # if lambda is not specified, the highest root is used
    if lamb == "none":
        lamb = lie_algebra.highest_root()
    else:
        lamb = changeBasis.inverse() * vector(lamb)
        lamb = weyl_group.domain()(list(eval(str(lamb))))
    
    # if mu is not specified, 0 weight is used
    if mu == "none":
        mu = weyl_group.domain()([0 for i in range(0, vectorDimension[name])])
    else:
        mu = changeBasis.inverse() * vector(mu)
        mu = weyl_group.domain()(list(eval(str(mu))))
    
    # rho is 1/2 the sum of the positive roots
    rho = lie_algebra.rho()
    maxl = 0
    for elm in weyl_group:
        length = elm.length()
        if (length > maxl):
            maxl = length
    
        if (length > 8):
            continue
        
        # compute the expression in the partition function
        res = elm.action(lamb + rho) - (mu + rho)
        
        # change basis from standard basis to simple roots
        res = vector(list(Eval(res)))
        res = changeBasis * res
        
        # make a weight object out of the result, check if it can be partitioned
        r = Weight(list(eval(str(res))))
        
        if not (r.isNegative() or r.hasFraction()):
            # root can be partitioned, find its Pq, add it to the sum
            tree = PartitionTree(r, alphaPoset, 0, 0)
            
            # delete unsucessful branches
            tree.cleanTree()
            
            # find value of Pq for this term
            # we want our vectors to be the same size, so we use the number of positive roots, and add extra for padding
            Pq = [0 for x in range (0, numOfRoots[name] + 25)]
            tree.generatePq(Pq)
            
            # save results to a table
            count.append([str(elm), str(elm.length()), str(r), str(Pq)])
            
            # determine the sign
            Pq = ((-1)**length) * vector(Pq)
            
            # add to the sum, for both multiplicity formulas
            equation = equation + Pq
            m += ((-1)**length) * tree.countPartitions()
            
            #delete the objects that we are no longer going to be using
            del(tree)
            del(Pq)
            
        del(res)
        
    return equation

# write all the collected information to a file for later parsing
f = open("E8 - results.txt", "w+")
for c in count:
    f.write(str(c) + "\n")
    
f.write("\n\nequation coeffs: " + str(computeMq("E8", "none", "none")) + "\n")
print(m)
f.close()
\end{lstlisting}

\subsection{The $q$-analog of Kostant's weight multiplicity formula by series expansion}
\begin{lstlisting}
# Anthony Simpson, als7@williams.edu
from sage.all import *
from sage.misc.parser import Parser
from Weight import Weight
import __future__

# we print to the file as we make our calculations
fl = open("results.txt", "w+")

# Here we define the variables that will appear in the expansion, since we are working with the exceptional Lie algebras we will need at most 8 A's and q, I used the first 8 letters of the alphabet
a, b, c, d, e, f, g, h, q = var('a, b, c, d, e, f, g, h, q')
# We build terms based on the simple roots, i.e. alpha1->A1, so we use this dictionary to easily make the correlation
variables = {'A1':a, 'A2':b, 'A3':c, 'A4':d, 'A5':e, 'A6':f, 'A7':g, 'A8':h, 'q':q}
P = Parser(make_var=var)

# Constants Declaration
bChangeA4 = matrix([[1, 0, 0, 0, 0],
                    [-1, 1, 0, 0, 0],
                    [0, -1, 1, 0, 0],
                    [0, 0, -1, 1, 0],
                    [0, 0, 0, -1, 1]]).inverse();

bChangeG2 = matrix([[2, 1, 0],
[1, 0, 0],
[1, 1, 1]])

bChangeF4 = matrix([[1, 1, 0, 0],
[2, 1, 1, 0],
[3, 1, 1, 1],
[2, 0, 0, 0]])

bChangeE6 = matrix([[0,   0,    0,    0,    0,   -2,    0,    0],
[ 1.0/2,  1.0/2,  1.0/2,  1.0/2,  1.0/2, -3.0/2,    0,    0],
[-1.0/2,  1.0/2,  1.0/2,  1.0/2,  1.0/2, -5.0/2,    0,    0],
[   0,    0,    1,    1,    1,   -3,    0,    0],
[   0,    0,    0,    1,    1,   -2,    0,    0],
[   0,    0,    0,    0,    1,   -1,    0,    0],
[   0,    0,    0,    0,    0,   -1,    1,    0],
[   0,    0,    0,    0,    0,    1,    0,    1]])

bChangeE7 = matrix([[   0,    0,    0,    0,    0,    0,   -2,    0],
[ 1.0/2,  1.0/2,  1.0/2,  1.0/2,  1.0/2,  1.0/2,   -2,    0],
[-1.0/2,  1.0/2,  1.0/2,  1.0/2,  1.0/2,  1.0/2,   -3,    0],
[   0,    0,    1,    1,    1,    1,   -4,    0],
[   0,    0,    0,    1,    1,    1,   -3,    0],
[   0,    0,    0,    0,    1,    1,   -2,    0],
[   0,    0,    0,    0,    0,    1,   -1,    0],
[   0,    0,    0,    0,    0,    0,    1,    1]])

bChangeE8 = matrix([[   0,    0,    0,    0,    0,    0,    0,    2],
[ 1.0/2,  1.0/2,  1.0/2,  1.0/2,  1.0/2,  1.0/2,  1.0/2,  5.0/2],
[-1.0/2,  1.0/2,  1.0/2,  1.0/2,  1.0/2,  1.0/2,  1.0/2,  7.0/2],
[   0,    0,    1,    1,    1,    1,    1,    5],
[   0,    0,    0,    1,    1,    1,    1,    4],
[   0,    0,    0,    0,    1,    1,    1,    3],
[   0,    0,    0,    0,    0,    1,    1,    2],
[   0,    0,    0,    0,    0,    0,    1,    1]])

basisChangeMatrices = {"G2": bChangeG2, "F4": bChangeF4, "E6": bChangeE6, "E7": bChangeE7, "E8": bChangeE8, "A4":bChangeA4}
numOfRoots = {"G2": 6, "F4": 4, "E6": 36, "E7": 63, "E8": 120}
vectorDimension = {"G2": 3, "F4": 4, "E6": 8, "E7": 8, "E8": 8}
num = 0

# this eval method is used instead of the built in eval function to fix some associated rounding errors
def Eval(string):
    return eval(compile(str(string), '<string>', 'eval', __future__.division.compiler_flag))

# this function takes a list of coefficients and transforms it into a geometic sum. For example, it takes the root a1 + 5a2 + 3a3 + 2a4 to (A1)(A2^5)(a3^3)(a4^2), then 1/(1 - (A1)(A2^5)(a3^3)(a4^2))
def toGeometricSum(weight):
    x = 1
    for i in range(0, len(weight)):
        for j in range(0, weight[i]):
            x = x * variables["A" + str(i+1)]
    return 1/(1 - q*x)

def computeMq(name, lamb, mu):
    # initialize constants and vector space for the lie algebra
    lie_algebra = RootSystem(name).ambient_space()
    weyl_group = WeylGroup(name, prefix = "s")

    # start mq = 0 then build it term by term
    equation = 0

    # used to change the basis from the standard basis of R^n to simple roots
    changeBasis = basisChangeMatrices[name]

    # this collect all of the poitive roots and builds the general terms for the summations
    poset = [vector(list(Eval(x))) for x in RootSystem(name).ambient_space().positive_roots()]
    alphaPoset = [changeBasis * x for x in poset]
    termsForSums = [toGeometricSum(list(x)) for x in alphaPoset]
    
    # create the geometric sum that we will be expanding by multiplying the terms created from the positive roots
    mainEq = 1
    for x in termsForSums:
        mainEq *= x

    # if lambda is not specified, the highest root is used
    if lamb == "none":
        lamb = lie_algebra.highest_root()
    else:
        lamb = changeBasis.inverse() * vector(lamb)
        lamb = weyl_group.domain()(list(eval(str(lamb))))

    # if mu is not specified, 0 vector is used
    if mu == "none":
        mu = weyl_group.domain()([0 for i in range(0, vectorDimension[name])])
    else:
        mu = changeBasis.inverse() * vector(mu)
        mu = weyl_group.domain()(list(eval(str(mu))))

    # rho is 1/2 the sum of the positive roots
    rho = lie_algebra.rho()

    for elm in weyl_group:
        # expression in partition function
        res = elm.action(lamb + rho) - (mu + rho)

        #change basis from standard basis to simple roots
        res = vector(list(Eval(res)))
        res = changeBasis * res
        r = Weight(list(res))

        if not (r.isNegative() or r.hasFraction()):
            print(str(res) + " - " + str(elm))

            answer = mainEq

            # expand the series until we reach one plus the power of the term that we are seeking.
            for i in range(0, len(res)):
                    answer = answer.series(variables["A" + str(i+1)], res[i] + 1).truncate()
                    # we take the coefficient of the desired power because that's all that matters, then expand that series to get the desired powers of the other variables
                    answer = answer.coefficient(variables["A" + str(i+1)], res[i])
            
            # at the end of the expanding, we will have the coefficient of the term corresponding to the weight we are working with, which will be exactly Pq of the weight we wanted
            answer = answer.expand()

            # print to a file for later use
            fl.write(str([str(elm), str(elm.length()), str(r), str(answer)]) + "\n")

            # add contribution to our running sum
            equation += answer * ((-1)**elm.length())

        del(res)
        del(r)

    return equation

#call the function and record the results
mq = computeMq("E8", "none", "none")


fl.write(str(mq))
fl.write("\n")
fl.close()

\end{lstlisting}

\subsection{The final program} \label{sect:fin}
\begin{lstlisting}
# Anthony Simpson, Williams College, als7@williams.edu
# this is used to stop rounding errors because of python2.7
def Eval(string):
    return eval(compile(str(string), '<string>', 'eval', __future__.division.compiler_flag))

def getBasisChange(name):
    # make a change of basis matrix by using the coordinate vectors of the simples in terms of the standard basis as row vectors
    simples = RootSystem(name).ambient_space().simple_roots()
    bChange = ([Eval((str(x).replace("(", "[").replace(")", "]"))) for x in simples])
    
    # may be a subspace, use the standard basis to fill in the matrix until square
    for i in range(len(bChange), len(bChange[0])):
        bChange.append([1 if j==i else 0 for j in range(0, len(bChange[0]))])
        
    # transpose so that simples are column vectors, basis change is now the inverse of this matrix
    return matrix(bChange).transpose().inverse()

def geometricSumForPartition(positive_root, translations, q_analog):
    # get the translation of the positive root to be used in a geometric sum
    x = 1
    for i in range(0, len(positive_root)):
        for j in range(0, positive_root[i]):
            # this gets the translation of each simple making the root
            x = x * translations["A" + str(i+1)]
            
    # return the correct sum based on if we are trying to compute the q-analog or not
    return 1/(1 - x) if not q_analog else 1/(1 -q*x)

# calculates the partition of a weight 
def calculatePartition(name, weight, positive_roots = [], translations = {}, q_analog = False):
    # if the set of positive roots is not given, we must get them
    if positive_roots == []:
        # get basis change matrix from standard to simples
        bChange = getBasisChange(name)
        
        # get the positive roots from Sage
        positive_roots = [vector(list(Eval(x))) for x in RootSystem(name).ambient_space().positive_roots()]
        
        # make those roots be in terms of the simple roots
        positive_roots = [bChange * x for x in positive_roots]
        
    # if the translations for the weights into variables are not provided, we make them here
    if translations == {}:
        s = 'A1'
        for i in range(1, len(weight)):
            s += ', A' + str(i + 1)
        variables = var(s)
        for i in range(0, len(weight)):
            translations["A" + str(i+1)] = eval("A" + str(i+1))
    
    # get the sum that will be expanded
    termsForSum = [geometricSumForPartition(list(x), translations, q_analog) for x in positive_roots]
    answer = 1
    for x in termsForSum:
        answer *= x
    
    # expand the sum in steps, only taking the portion that will be useful
    for i in range(0, len(weight)):
        answer = answer.series(translations["A" + str(i+1)], weight[i] + 1).truncate()
        answer = answer.coefficient(translations["A" + str(i+1)], weight[i])
    
    # return final answer
    answer = answer.expand()
    return answer
    
def findAltSet(name, lamb = None, mu = None):
    # initialize constants and vector space for the lie algebra
    lie_algebra = RootSystem(name).ambient_space()
    weyl_group = WeylGroup(name, prefix = "s")
    simples = weyl_group.gens()
        
    altset = [weyl_group.one()]

    # used to change the basis from the standard basis of R^n to simple roots
    changeBasis = getBasisChange(name)

    # if lambda is not specified, the highest root is used
    if lamb == None:
        lamb = lie_algebra.highest_root()
    else:
        lamb = changeBasis.inverse() * vector(lamb)
        lamb = weyl_group.domain()(list(eval(str(lamb))))

    # if mu is not specified, 0 vector is used
    if mu == None:
        mu = weyl_group.domain()([0 for i in range(0, len(lie_algebra.simple_roots()))])
    else:
        mu = changeBasis.inverse() * vector(mu)
        mu = weyl_group.domain()(list(eval(str(mu))))

    rho = lie_algebra.rho()
    length = len(altset)
    i=0
    # this loops over our running subset of the alternation set. We start the alternation set with just the identity in it, then we build it by appending simple reflections, until that applied to our expression results in a weight that can be partitioned
    while i < length:
        for simple in simples:
            # don't want to add things that will result in repeats (i.e. s1*s1 = 1) as it will make us loop infinitely
            if ((altset[i] == simple) or (altset[i] == altset[i] * simple)):
                continue
                
            # compute the result of applying the reflection to the expressino
            res = (altset[i]*simple).action(lamb + rho) - (rho + mu)
            res = changeBasis * vector(list(Eval(res)))
            res = Weight(res)
            
            # if the result can be partitioned, add the reflection to the running alternation set
            if not (res.isNegative() or res.hasFraction()):
                if not (altset[i]*simple in altset):
                    altset.append(altset[i]*simple)
                    length += 1
        i+=1
        
    return altset
    

def calculateMultiplicity(name, lamb = None, mu = None, q_analog = False):
    # initializing variables that will be used in this function
    mult = 0
    lie_algebra = RootSystem(name).ambient_space()
    weyl_group = WeylGroup(name, prefix = "s")
    
    # used to change the basis from the standard basis of R^n to simple roots
    changeBasis = getBasisChange(name)
    
    # get positive roots in terms of the simples
    positive_roots = [vector(list(Eval(x))) for x in RootSystem(name).ambient_space().positive_roots()]
    positive_roots = [getBasisChange(name) * x for x in positive_roots]

    # if lambda is not specified, the highest root is used
    if lamb == None:
        lamb = lie_algebra.highest_root()
    else:
        lamb = changeBasis.inverse() * vector(lamb)
        lamb = weyl_group.domain()(list(eval(str(lamb))))

    # if mu is not specified, 0 vector is used
    if mu == None:
        mu = weyl_group.domain()([0 for i in range(0, len(lie_algebra.simple_roots()))])
    else:
        mu = changeBasis.inverse() * vector(mu)
        mu = weyl_group.domain()(list(eval(str(mu))))

    rho = lie_algebra.rho()
    
    # we know we only have to iterate over the elements of the alternation set, so we find it and begin
    altset = findAltSet(name)
    for elm in altset:
        # expression in partition function
        res = elm.action(lamb + rho) - (mu + rho)

        # change basis from standard basis to simple roots
        res = vector(list(Eval(res)))
        res = changeBasis * res
        
        # calculate the partiton, or the q-analog, based on the settings
        term = calculatePartition(name, list(res), positive_roots, q_analog=q_analog)
        
        # calcuate the contribution to the sum and add it
        term *= (-1)**elm.length()
        mult += term
        
    return mult
\end{lstlisting}

\newpage
\section{Application of program with associated tables}\label{sec:application}

\subsection{Alternation set for $G_2$}

\begin{longtable}{|c|c|c|p{4cm}|p{6cm}|}
\rowcolor{lightgray}
\rowcolors{1}{}{midgray}
No. \ &$\sigma$ \ &$\ell(\sigma)$ \ &$\xi = \sigma(\widetilde{\alpha} + \rho) - \rho$ \ &$\wp_q(\xi)$\\

1 & $1$ & 0 & $ 3\alpha_{1} + 2\alpha_{2} $ & $ q^{1} + 2q^{2} + 2q^{3} + q^{4} + q^{5} $\\\hline
2 & $s_1$ & 1 & $ 2\alpha_{1} + 2\alpha_{2} $ & $ 2q^{2} + q^{3} + q^{4} $\\\hline
3 & $s_2$ & 1 & $ 3\alpha_{1} $ & $ q^{3} $\\\hline

\rowcolor{white}\multicolumn{5}{|c|}{$m_q(\widetilde{\alpha}, 0) = q + q^5$}\\

\hline
\end{longtable}

\subsection{Alternation set for $F_4$}

\begin{longtable}{|c|c|c|p{4cm}|p{6cm}|}
\rowcolor{lightgray}
\rowcolors{1}{}{midgray}
No. \ &$\sigma$ \ &$\ell(\sigma)$ \ &$\xi = \sigma(\widetilde{\alpha} + \rho) - \rho$ \ &$\wp_q(\xi)$ \\
1 & $1$ & 0 & $ 2\alpha_{1} + 3\alpha_{2} + 4\alpha_{3} + 2\alpha_{4} $ & $ q^{1} + 7q^{2} + 27q^{3} + 53q^{4} + 69q^{5} + 59q^{6} + 40q^{7} + 20q^{8} + 9q^{9} + 3q^{10} + q^{11} $\\\hline
2 & $s_1$ & 1 & $ 3\alpha_{2} + 4\alpha_{3} + 2\alpha_{4} $ & $ 5q^{3} + 9q^{4} + 13q^{5} + 8q^{6} + 5q^{7} + 2q^{8} + q^{9} $\\\hline
3 & $s_2$ & 1 & $ 2\alpha_{1} + 2\alpha_{2} + 4\alpha_{3} + 2\alpha_{4} $ & $ 3q^{2} + 10q^{3} + 26q^{4} + 31q^{5} + 29q^{6} + 17q^{7} + 9q^{8} + 3q^{9} + q^{10} $\\\hline
4 & $s_3$ & 1 & $ 2\alpha_{1} + 3\alpha_{2} + 3\alpha_{3} + 2\alpha_{4} $ & $ 3q^{2} + 16q^{3} + 33q^{4} + 41q^{5} + 34q^{6} + 20q^{7} + 9q^{8} + 3q^{9} + q^{10} $\\\hline
5 & $s_4$ & 1 & $ 2\alpha_{1} + 3\alpha_{2} + 4\alpha_{3} + \alpha_{4} $ & $ 3q^{2} + 16q^{3} + 32q^{4} + 36q^{5} + 29q^{6} + 17q^{7} + 8q^{8} + 3q^{9} + q^{10} $\\\hline
6 & $s_2s_1$ & 2 & $ 4\alpha_{3} + 2\alpha_{4} $ & $ q^{4} + q^{5} + q^{6} $\\\hline
7 & $s_2s_3$ & 2 & $ 2\alpha_{1} + \alpha_{2} + 3\alpha_{3} + 2\alpha_{4} $ & $ 2q^{3} + 6q^{4} + 8q^{5} + 6q^{6} + 3q^{7} + q^{8} $\\
\hline
\end{longtable}

\subsection{Alternation set for $E_6$}

\begin{longtable}[H]{|c|c|c|p{3.5cm}|p{6cm}|}
\rowcolor{lightgray}
\rowcolors{1}{}{midgray}
No. \ &$\sigma$ \ &$\ell(\sigma)$ \ &$\xi = \sigma(\widetilde{\alpha} + \rho) - \rho$ \ &$\wp_q(\xi)$ \\
1 & $1$ & 0 & $ \alpha_{1} + 2\alpha_{2} + 2\alpha_{3} + 3\alpha_{4} + 2\alpha_{5} + \alpha_{6} $ & $ q^{1} + 10q^{2} + 45q^{3} + 105q^{4} + 150q^{5} + 142q^{6} + 97q^{7} + 48q^{8} + 18q^{9} + 5q^{10} + q^{11} $\\\hline
2 & $s_1$ & 1 & $ 2\alpha_{2} + 2\alpha_{3} + 3\alpha_{4} + 2\alpha_{5} + \alpha_{6} $ & $ 3q^{2} + 18q^{3} + 43q^{4} + 57q^{5} + 50q^{6} + 30q^{7} + 13q^{8} + 4q^{9} + q^{10} $\\\hline
3 & $s_2$ & 1 & $ \alpha_{1} + 2\alpha_{3} + 3\alpha_{4} + 2\alpha_{5} + \alpha_{6} $ & $ 6q^{3} + 18q^{4} + 26q^{5} + 20q^{6} + 11q^{7} + 4q^{8} + q^{9} $\\\hline
4 & $s_3$ & 1 & $ \alpha_{1} + 2\alpha_{2} + \alpha_{3} + 3\alpha_{4} + 2\alpha_{5} + \alpha_{6} $ & $ 3q^{2} + 18q^{3} + 46q^{4} + 66q^{5} + 60q^{6} + 37q^{7} + 16q^{8} + 5q^{9} + q^{10} $\\\hline
5 & $s_4$ & 1 & $ \alpha_{1} + 2\alpha_{2} + 2\alpha_{3} + 2\alpha_{4} + 2\alpha_{5} + \alpha_{6} $ & $ 3q^{2} + 18q^{3} + 47q^{4} + 69q^{5} + 65q^{6} + 41q^{7} + 18q^{8} + 5q^{9} + q^{10} $\\\hline
6 & $s_5$ & 1 & $ \alpha_{1} + 2\alpha_{2} + 2\alpha_{3} + 3\alpha_{4} + \alpha_{5} + \alpha_{6} $ & $ 3q^{2} + 18q^{3} + 46q^{4} + 66q^{5} + 60q^{6} + 37q^{7} + 16q^{8} + 5q^{9} + q^{10} $\\\hline
7 & $s_6$ & 1 & $ \alpha_{1} + 2\alpha_{2} + 2\alpha_{3} + 3\alpha_{4} + 2\alpha_{5} $ & $ 3q^{2} + 18q^{3} + 43q^{4} + 57q^{5} + 50q^{6} + 30q^{7} + 13q^{8} + 4q^{9} + q^{10} $\\
\hline
\end{longtable}

\subsection{Alternation set for $E_7$}

\begin{longtable}[H]{|c|c|c|p{3cm}|p{6cm}|}
\rowcolor{lightgray}
\rowcolors{1}{}{midgray}
No. \ &$\sigma$ \ &$\ell(\sigma)$ \ &$\xi = \sigma(\widetilde{\alpha} + \rho) - \rho$ \ &$\wp_q(\xi)$ \\
1 & $1$ & 0 & $ 2\alpha_{1} + 2\alpha_{2} + 3\alpha_{3} + 4\alpha_{4} + 3\alpha_{5} + 2\alpha_{6} + \alpha_{7} $ & $ q^{1} + 16q^{2} + 120q^{3} + 500q^{4} + 1340q^{5} + 2488q^{6} + 3420q^{7} + 3617q^{8} + 3054q^{9} + 2106q^{10} + 1216q^{11} + 592q^{12} + 247q^{13} + 87q^{14} + 26q^{15} + 6q^{16} + q^{17} $\\\hline
2 & $s_1$ & 1 & $ 2\alpha_{2} + 3\alpha_{3} + 4\alpha_{4} + 3\alpha_{5} + 2\alpha_{6} + \alpha_{7} $ & $ 15q^{3} + 90q^{4} + 275q^{5} + 500q^{6} + 642q^{7} + 610q^{8} + 458q^{9} + 273q^{10} + 136q^{11} + 55q^{12} + 19q^{13} + 5q^{14} + q^{15} $\\\hline
3 & $s_2$ & 1 & $ 2\alpha_{1} + \alpha_{2} + 3\alpha_{3} + 4\alpha_{4} + 3\alpha_{5} + 2\alpha_{6} + \alpha_{7} $ & $ 4q^{2} + 44q^{3} + 218q^{4} + 628q^{5} + 1198q^{6} + 1630q^{7} + 1674q^{8} + 1346q^{9} + 877q^{10} + 471q^{11} + 213q^{12} + 80q^{13} + 25q^{14} + 6q^{15} + q^{16} $\\\hline
4 & $s_3$ & 1 & $ 2\alpha_{1} + 2\alpha_{2} + 2\alpha_{3} + 4\alpha_{4} + 3\alpha_{5} + 2\alpha_{6} + \alpha_{7} $ & $ 4q^{2} + 44q^{3} + 218q^{4} + 628q^{5} + 1216q^{6} + 1684q^{7} + 1764q^{8} + 1437q^{9} + 944q^{10} + 506q^{11} + 227q^{12} + 84q^{13} + 26q^{14} + 6q^{15} + q^{16} $\\\hline
5 & $s_4$ & 1 & $ 2\alpha_{1} + 2\alpha_{2} + 3\alpha_{3} + 3\alpha_{4} + 3\alpha_{5} + 2\alpha_{6} + \alpha_{7} $ & $ 4q^{2} + 44q^{3} + 218q^{4} + 634q^{5} + 1228q^{6} + 1709q^{7} + 1797q^{8} + 1480q^{9} + 980q^{10} + 531q^{11} + 237q^{12} + 87q^{13} + 26q^{14} + 6q^{15} + q^{16} $\\\hline
6 & $s_5$ & 1 & $ 2\alpha_{1} + 2\alpha_{2} + 3\alpha_{3} + 4\alpha_{4} + 2\alpha_{5} + 2\alpha_{6} + \alpha_{7} $ & $ 4q^{2} + 44q^{3} + 216q^{4} + 624q^{5} + 1201q^{6} + 1658q^{7} + 1731q^{8} + 1413q^{9} + 931q^{10} + 502q^{11} + 226q^{12} + 84q^{13} + 26q^{14} + 6q^{15} + q^{16} $\\
\hline
\end{longtable}

\newpage
\subsection{Alternation set for $E_8$}

\begin{longtable}[H]{|c|c|c|p{3cm}|p{5.1cm}|}
\rowcolor{lightgray}
\rowcolors{1}{}{midgray}
No. \ &$\sigma$ \ &$\ell(\sigma)$ \ &$\xi = \sigma(\widetilde{\alpha} + \rho) - \rho$ \ &$\wp_q(\xi)$ \\
1 & $1$ & 0 & $ 2\alpha_{1} + 3\alpha_{2} + 4\alpha_{3} + 6\alpha_{4} + 5\alpha_{5} + 4\alpha_{6} + 3\alpha_{7} + 2\alpha_{8} $ & $ q+ 28q^{2 }+ 378q^{3 }+ 2961q^{4 }+ 15435q^{5 }+ 57540q^{6 }+ 161931q^{7 }+ 357494q^{8 }+ 638680q^{9 }+ 946515q^{10 }+ 1188475q^{11 }+ 1287007q^{12 }+ 1220804q^{13 }+ 1027773q^{14 }+ 776864q^{15 }+ 532288q^{16 }+ 333340q^{17 }+ 192015q^{18 }+ 102283q^{19 }+ 50542q^{20 }+ 23216q^{21 }+ 9902q^{22 }+ 3913q^{23 }+ 1421q^{24 }+ 470q^{25 }+ 138q^{26 }+ 35q^{27 }+ 7q^{28 }+q^{29 }$\\\hline
2 & $s_1$ & 1 & $ \alpha_{1} + 3\alpha_{2} + 4\alpha_{3} + 6\alpha_{4} + 5\alpha_{5} + 4\alpha_{6} + 3\alpha_{7} + 2\alpha_{8} $ & $ 6q^{2}+ 126q^{3 }+ 1220q^{4 }+ 7125q^{5 }+ 28306q^{6 }+ 82264q^{7 }+ 183855q^{8 }+ 327849q^{9 }+ 480231q^{10 }+ 591806q^{11 }+ 626034q^{12 }+ 578168q^{13 }+ 472913q^{14 }+ 346746q^{15 }+ 230214q^{16 }+ 139536q^{17 }+ 77721q^{18 }+ 39965q^{19 }+ 19026q^{20 }+ 8387q^{21 }+ 3416q^{22 }+ 1277q^{23 }+ 434q^{24 }+ 131q^{25 }+ 34q^{26 }+ 7q^{27 }+q^{28 }$\\\hline
3 & $s_2$ & 1 & $ 2\alpha_{1} + 2\alpha_{2} + 4\alpha_{3} + 6\alpha_{4} + 5\alpha_{5} + 4\alpha_{6} + 3\alpha_{7} + 2\alpha_{8} $ & $ 6q^{2}+ 126q^{3 }+ 1223q^{4 }+ 7188q^{5 }+ 28838q^{6 }+ 84755q^{7 }+ 191690q^{8 }+ 345879q^{9 }+ 512388q^{10 }+ 637867q^{11 }+ 680643q^{12 }+ 632993q^{13 }+ 520542q^{14 }+ 383122q^{15 }+ 255012q^{16 }+ 154760q^{17 }+ 86204q^{18 }+ 44257q^{19 }+ 21003q^{20 }+ 9208q^{21 }+ 3723q^{22 }+ 1377q^{23 }+ 462q^{24 }+ 137q^{25 }+ 35q^{26 }+ 7q^{27 }+q^{28 }$\\
\hline
\end{longtable}

\end{document}